\documentclass[12pt,titlepage]{article}
\usepackage{amsmath,amsthm, amssymb}
\usepackage{multicol}
\usepackage{color}
\usepackage{secdot}
\allowdisplaybreaks[4]
\makeatletter
\def\section{\@startsection{section}{1}{\z@}{-4ex plus -1ex minus -.2ex}{1.5ex plus .2ex}{\Large\bf}}
\makeatother
\makeatletter

\@addtoreset{equation}{section}     
\makeatother
\def\because{\raisebox{1.1ex}{.}.\raisebox{1.1ex}{.}}

\def\R{\mathbf{R}}

\newtheorem{lemma}{Lemma}
\newtheorem{proposition}{Proposition}
\newtheorem{theorem}{Theorem}

\theoremstyle{definition}

\newtheorem{remark}{Remark}

\makeatletter
\def\@proofcounterend{.}
\def\@proofcounterend{}
\def\proof{\@ifnextchar[{\@yproof}{\@xproof}}
\def\@xproof{\@beginproof\ignorespaces}
\def\@yproof[#1]{\@opargbeginproof{#1}\ignorespaces}
\def\@beginproof{\pushQED{\qed}
\topsep6\p@\@plus6\p@\relax
\begin{trivlist}
  \item[\hskip\labelsep{\bf ???\@proofcounterend}]\rm}
\def\@opargbeginproof#1{\pushQED{\qed}
\topsep6\p@\@plus6\p@\relax
\begin{trivlist}
  \item[\hskip\labelsep{\bf #1????\@proofcounterend}]\rm}
\def\endproof{\popQED \end{trivlist}}  
\def\subproof{\@beginsubproof\ignorespaces}
\def\beginsubproof{
 \topsep6\p@\@plus6\p@\relax
\topsep4\p@\@plus4\p@\relax
\begin{trivlist}
  \item[\hskip\labelsep{
  \mbox{\ooalign
  {\hfill {\scriptsize $\because$}\hfill\crcr $\bigcirc$}}}]\rm}
\def\endsubproof{\end{trivlist}}
\makeatother                    
\makeatletter
\@addtoreset{claim}{section}
\@addtoreset{definition}{section}
\@addtoreset{theorem}{section}
\@addtoreset{lemma}{section}
\@addtoreset{remark}{section}
\@addtoreset{example}{section}
\@addtoreset{corollary}{section}
\@addtoreset{proposition}{section}
\makeatletter
\newcommand{\Equref}[1]{(\ref{#1})}  
\makeatother
\begin{document}

\begin{center} \large A Probabilistic Approach to the Zero--Mass Limit Problem \\for Three Magnetic Relativistic Schr\"odinger Heat Semigroups 
\end{center}

\begin{center}
Taro  Murayama
\end{center}

Abstract: We consider three    magnetic relativistic Schr\"odinger operators which correspond to  the same classical symbol $\sqrt{(\xi-A(x))^2+m^2}+V(x)$ and whose heat semigroups admit   the  Feynman-Kac-It\^o type    path integral  representation $E[e^{-S^m(x,t; X)}g(x+X(t))]$.  Using   these representations, we  prove  the convergence of these  heat  semigroups %by using their path integral representations 
  when    the     mass--parameter $m$    goes to zero. Its  proof  reduces to   the convergence  of          $e^{-S^m(x,t;X)}$, which  yields  a limit theorem for   exponentials of   semimartingales as functionals of L\'evy processes $X$.

\section{Introduction and Results} 
In  a recent paper \cite{I and M 2014},  we studied  the {\it zero--mass limit problem}  for heat semigroup of the   Weyl--pseudodifferential  operator $H_{A}^m+V$ with classical symbol $\sqrt{(\xi-A(x))^2+m^2}+V(x)$ to show that as $m\downarrow 0$, 
\begin{align}
e^{-t[H_{A}^m-m+V]}\to e^{-t[H_A^0+V]}\quad \text{\rm strongly},
\end{align}
uniformly on every finite bounded  interval in $t\geq0$.   For the proof,  its {\it Feynman--Kac--It\^o (F--K--I)}   type path integral formula (e.g. \cite{S 79}) was used. 
Here $m$ is the {\it mass parameter}, and
$A:\R^d \rightarrow \R^d, \, V:\R^d \rightarrow \R$ are the magnetic vector  potential,
the electric scalar potential, which in fact were assumed to satisfy that
$A \in  C_0^{\infty}(\R^d;\R^d)$, $V \in C_0(\R^d;\R)$.

In this paper,  we study this problem under more general assumptions on the potentials
$A$ and $V$, and treat moreover the additional case for
other two different magnetic relativistic Schr\"odinger operators $H_{A}^m+V$,
together with their respective F--K--I   type    formulae.
The problem will be solved by  discussing  the convergence of special kind of semimartingales,   namely, {\it exponentials of semimartingales}, as functionals of     {\it L\'evy processes} (see Lemma 5.1 and Lemma 6.1).  To best my knowledge, such convergence  does not seem to  have been treated in the framework of the limit theorems for semimartingales represented by stochastic integrals (cf. \cite{K and W 86}).

Now, let  $H_{1,A}^m$, $H_{2,A}^m$, $H_{3,A}^m$   denote  the following three magnetic relativistic Schr\"odinger operators   corresponding to the symbol of the classical kinetic energy    $\sqrt{(\xi-A(x))^2+m^2}$ ($(\xi, x)\in \R^d\times \R^d$):   
\begin{align}
  (H^{m}_{1, A}f)(x)    &:=\frac{1}{(2\pi)^{d} }    \iint_{\R^{d}\times \R^{d}} 
e^{i(x-y)\cdot \xi}\sqrt{\left(\xi -A(\tfrac{x+y}{2})\right)^{2}+m^{2}}  f(y)dy d\xi,\\            %\quad f \in  C_{0}^{\infty}(\mathbf{R}^{d}). 
 (H^{m}_{2, A}f)(x)    &:=\frac{1}{(2\pi)^{d} }    \iint_{\R^{d}\times \R^{d}} 
e^{i(x-y)\cdot \xi}\sqrt{\left(\xi -\int_0^1A((1-\theta)x+\theta y )d\theta      \right)^{2}+m^{2}}  f(y)dy d\xi, \\%\quad f \in  C_{0}^{\infty}(\mathbf{R}^{d}). 
  H_{3,A}^m&: = \sqrt{(-i\nabla-A(x))^2+m^2}. 
\end{align}
$H_{1,A}^m$ is the {\it  Weyl pseudodifferential operator  
} introduced  in      \cite{I and T 86} and 
studied further  in   \cite{I 2012}, \cite{I 2013}. $H_{2,A}^m$ is the {\it
pseudodifferential operator}      
  defined as a modification of $H_{1,A}^m$     (\cite{I and M and P 2007}, \cite{I and M and P 2008}, \cite{I and M and P 2010}).    $H_{3,A}^m$ is {\it the square root of  the       nonnegative selfadjoint operator}  $(-i\nabla-A(x))^2+m^2$ in $L^2(\R^d)$. Each operator $H_{j,A}^m+V$ ($j=1,2,3$) may be used to        describe the motion of    a  relativistic spinless particle with  mass $m \geq 0$  in  the electromagnetic field.  We have  $H_{1,0}^m=H_{2,0}^m=H_{3,0}^m=\sqrt{-\Delta+m^2}$ for $A \equiv 0$, where $\Delta$ is the Laplacian in  $\R^d$.  For $A \not \equiv 0 $, the operators 
$H_{1,A}^m$, $H_{2,A}^m$, $H_{3,A}^m$ are different   from one another, although they coincide in the case of constant magnetic field, i.e.,  when $A(x)=\dot{A}x$ with $\dot{A}$ a constant symmetric   matrix.    
Under gauge  transformation,        $H_{2,A}^m$ and $H_{3,A}^m$ are 
covariant, but $H_{1,A}^m$ is not (\cite[Section 2, Section 3]{I 2012}, \cite[Section 2]{I 2013}).      

Let us consider the  heat semigroups  $e^{-t[H_{j, A}^m-m+V]}g$ applied to a function $g$,  each of  which is      the solution $u(x,t)= ( e^{-t[H_{j, A}^m-m+V]}g)(x)$       of the Cauchy problem for the heat equation 
\begin{align*}
\begin{cases}
& \frac{\partial}{\partial t}u(x,t)=-[H_{j, A}^m-m+V]u(x,t),\quad x\in \R^d,\ t>0,\\
& u(x,0)=g(x),\quad x\in \R^d.
\end{cases}
\end{align*}
They are  known (\cite{I and T 86}, \cite{I and M and P 2007}, \cite{ D and S 90})   to   be represented by F-K-I    type formulae  as follows: 
\begin{align}
(e^{-t[H_{j,A}^m-m+V]}g)(x)& = E^{\lambda^m}\bigg[ e^{-S^{m}_{j,A,V}(x,t; X)}g(x+X(t)) \bigg], \quad j=1,2, \\% d\lambda^m(X),\\
% (e^{-t[H_{2,A}^m-m+V]}g)(x)&=  E^{\lambda^m}\bigg[e^{-S^{m}_{2,A,V}(x,t; X)}g(x+X(t))\bigg],\\% d\lambda^m(X),\\
  (e^{-t[H_{3,A}^m-m+V]}g)(x)&  = E^{\mu\times \nu^m}\bigg[e^{-S_{3,A,V}       (x,  t; B, T)}g\Big(x+B\big(T(t)\big)\Big) \bigg].
\end{align}
Here we denote    by $E^{\mathbf{P}}[\cdots]=\int \cdots  d\mathbf{P}$ the expectation with respect to    the probability    measure $\mathbf{P}$. 
$\lambda^m$ and $\nu^m$ are some probability measures connected with $d$-dimensional  L\'{e}vy process $X$ and 1-dimentional subordinator $T$ to be   introduced as {\it time change}, respectively.  $\mu$ is the $d$-dimensional Wiener measure  associated with      $d$-dimensional standard Brownian motion $B$. 
 $S_{1,A,V}^m(x,t;X)$, $S_{2,A,V}^m(x,t;X)$   and   $S_{3,A,V} (x,t;B,T)$ are   complex-valued  {\it semimartingales}  given by stochastic  integrals of potentials  $A$ and $V$.

\vspace{4mm}

Our first result is  the  weak convergence of  two    probability measures $\lambda^m$ and $\nu^m$. 
\begin{theorem}
{\rm (i)} 
$\lambda^m$  weakly  converges  to $\lambda^0$ as $m \downarrow 0$.\\
{\rm (ii)} $\nu^m$  weakly converges to $\nu^0$ as $m \downarrow 0$.  
\end{theorem}
Our second result is the   strong convergence   of the   heat semigroups $e^{-t[H_{j,A}^m-m+V]}$ ($j=1,2,3$)  on   $C_\infty(\R^d):=\Big\{g \in C(\R^d); \displaystyle \lim_{|x|\to \infty}g(x)=0\Big\}$ with norm $\|g\|_\infty:=\displaystyle \sup_{x\in\R^d}|g(x)|$: 
\begin{align}\label{convergence in sup norm}
\sup_{t\leq t_0}\| e^{-t[H_{j,A}^m-m+V]}g-e^{-t[H_{j,A}^0+V]}g\|_\infty \to 0 \quad \text{as }\ m \downarrow 0. 
\end{align}

\begin{theorem} Assume that $g\in C_\infty(\R^d)$ and  $0\leq V \in C(\R^d;\R)$. \\ 
 {\rm(i)} 
If $A$ is locally $\alpha$-H\"older continuous  $(0<\alpha \leq 1)$, then \Equref{convergence in sup norm} holds  for $j=1,2$. 
\\
{\rm(ii)} If $A \in C^1(\R^d;\R^d)$, then       \Equref{convergence in sup norm}  holds                                for $j=3$. 
\end{theorem}

Our third result is the  strong convergence   of the  heat    semigroups $e^{-t[H_{j,A}^m-m+V]}$ ($j=1,2,3$)  on      $L^2(\R^d)$: 
\begin{align}\label{convergence in ell 2 norm}
\sup_{t\leq t_0}\| e^{-t[H_{j,A}^m-m+V]}g-e^{-t[H_{j,A}^0+V]}g\|_2 \to 0 \quad \text{as }\ m \downarrow 0. 
\end{align}
\begin{theorem} Assume that $g\in L^2(\R^d)$ and  $0\leq V \in L^1_\text{loc}(\R^d;\R)$. \\ 
 {\rm(i)} 
If $A\in    L^{1+\delta}_{\text{loc}}(\R^d;\R^d)$ for some $\delta>0$,  then   \Equref{convergence in ell 2 norm} holds for $j=1,2$. 
\\
{\rm(ii)} If $A\in  L^2_{\text{loc}}(\R^d;\R^d)$ and $\nabla \cdot A \in L^1_{\text{loc}}(\R^d;\R)$, then  \Equref{convergence in ell 2  norm}  holds for $j=3$. 
\end{theorem}

Claim (i) for $j=1$ of Theorem 1.2 and Theorem 1.3    are direct  generalizations of those results of \cite{I and M 2014}.  
Note that Theorem 1.2 and Theorem 1.3 hold if $V$ is {\it bounded from below}.   In fact, we have  only to replace $V$ by $V-\inf V$($\geq 0$). 

The problem may be thought of  for the operators $H_{j,A}^m +V$,  which are bounded from
below with more general  scalar potential  $V(x)$, for instance, a negative Coulomb potential $V(x) = - \frac{c}{|x|}$.
In fact, $H_{j,0}^m - \frac{c}{|x|}=\sqrt{-\Delta+m^2} - \frac{c}{|x|}$, with $c \leq \frac{2}{\pi}$, is known to be, as a quadratic form, bounded from below (nonnegative).
However, in this paper we content ourselves only with treating the above mentioned case,
partly because the mass parameter $m$ is involved only with the kinetic energy part $H_{j,A}^m$ containing
vector potential $A$ but not with scalar potential $V$, and partly because of
avoiding inessential difficulty coming from negativity of scalar potential.

\bigskip
This paper  is organized as follows.  In Section 2,  we  will describe  more precisely  
 the three F-K-I type formulae   (1.5) and  (1.6). 
In Section 3, we prove Theorem 1.1. In Section 4, we give preliminaries     to prove Theorem 1.2 and Theorem 1.3.   
In Section 5 and Section 6, we prove Theorem 1.2 and  Theorem 1.3, respectively.

\section{Three F--K--I type formulae} 
In this section, we give  more precise description of   the  three   F-K-I type   formulae  (1.5) and (1.6).   

For (1.5), 
$\lambda^m$ is the probability measure on the path space 
\begin{align*}
D_0=D_0([0,\infty)\to \R^d):= \{X:[0,\infty)\to \R^d;  X \text{ is c\`adl\`ag}, X(0)=0\}, 
\end{align*}
and satisfies  
\begin{align}\label{characteristic f.t. of lambda}
E^{\lambda^m}\big[  e^{i\xi \cdot X(t)}\big]  =e^{-t[\sqrt{\xi^2+m^2}-m]},\quad \xi \in \R^d, \ t\geq 0. 
\end{align}
 $X$ is a pure--jump L\'evy process with respect to $\lambda^m$, i.e.,                       
\begin{align}\label{eq:eq levy-ito}
X(t) =\int_{0}^{t}  \int_{|y| \geq 1}  y \ N_{X}(dsdy) 
+ \int_{0}^{t} \int_{0< |y| < 1}  y\ \widetilde{N_{X}^{m}}(dsdy),\quad \lambda^{m}\text{-a.s}. %X \in D_{0}. 
\end{align}
 $N_X(dsdy)$ is the counting measure on $(0, \infty) \times  \{|y|>0\}$ defined by 
\begin{align*}
N_X(dsdy):= \#  \{ u \in ds;  X(u)-X(u-)\in  dy\}. 
\end{align*} 
It  is the stationary  Poisson random measure with intensity measure (compensator)  $dsn^m(dy)$ with respect to $\lambda^m$, where     $n^m(dy)$ is the L\'evy  measure having   density 
\begin{align}\label{eq:eq levy measure}
n^{m}(y)
=  \begin{cases}
\displaystyle    2\left (\frac{m}{2\pi}\right)^{\frac{d+1}{2}}\frac{K_{\frac{d+1}{2}}(m|y|)}{|y|^{\frac{d+1}{2}}}, 
\quad & m>0, 
\\
 \displaystyle      \frac{\Gamma (\tfrac{d+1}{2})}{\pi^{\frac{d+1}{2}}} \frac{1}{|y|^{d+1}}, & m=0,
\end{cases}
\end{align}
so that $n^m(dy)=n^m(y)dy$, 
and then satisfies
\begin{align*}  
\int_{0<|y|<1}|y|^{1+\delta} n^m(dy)<\infty, \quad \delta>0,\ m\geq 0.
\end{align*} 
Here     $K_{\nu}$ stands  for the modified Bessel function of the third kind of order $\nu$ and $\Gamma$ denotes  the gamma function, respectively. 
$\widetilde{N_X^m}(dsdy)$ is the compensated   Poisson random measure, i.e.,  
\begin{align*}
\widetilde{N_X^m}(dsdy):=N_X(dsdy)-dsn^m(dy).
\end{align*}
$S^{m}_{1,A,V} (x,t; X)$ is a complex-valued semimartingale given    by  
 \begin{align}\label{semi 1}
S^{m}_{1,A,V} (x,t; X)&  :=   i  \bigg[  \int_{0}^{t}\int_{|y| \geq 1} A(x+X(s-)+\tfrac{1}{2}y) \cdot y N_{X}(dsdy) \nonumber \\
& \quad \quad + \int_{0}^{t}  \int_{0<|y| < 1}   A(x+X(s-)+\tfrac{1}{2}y) \cdot y \widetilde{N_{X}^{m}}(dsdy)\nonumber \\
 & \quad\quad + \int_{0}^{t} ds \text{ p.v.} \int_{0<|y|<1}  A(x+X(s)+\tfrac{1}{2}y) \cdot y n^{m}(dy)  \bigg]  \nonumber   \\
& \quad \quad + \int_{0}^{t} V(x+X(s))ds.
\end{align}
Here \lq\lq p.v.\rq\rq \   means        the principal value integral.   $S^{m}_{2,A,V}(x,t; X)$ is given           by  a modification of  $S^m_{1,A,V}(x,t;X)$
as follows:
\begin{align}\label{semi 2}
S^{m}_{2,A,V} (x,t; X) 
& : =   i  \bigg[  \int_{0}^{t}\int_{|y| \geq 1} \left(\int_0^1A(x+X(s-)+\theta y) d\theta\right) \cdot y N_{X}(dsdy) \nonumber \\
 & \quad\quad + \int_{0}^{t}  \int_{0<|y| < 1}   \left(\int_0^1A(x+X(s-)+\theta y) d\theta\right) \cdot y \widetilde{N_{X}^{m}}(dsdy)\nonumber \\
 & \quad\quad + \int_{0}^{t} ds \text{ p.v.} \int_{0<|y|<1} \left(\int_0^1A(x+X(s)+\theta y) d\theta\right) \cdot y n^{m}(dy)  \bigg]  \nonumber   \\
& \quad\quad  + \int_{0}^{t} V(x+X(s))ds. 
\end{align}

For (1.6), 
$\mu$ is the $d$-dimensional Wiener measure associated with $d$-dimensional   standard   Brownian  motion $B$. $\nu^m$ is the probablity measure    on 
\begin{align*}
D_0^{(1)}&=D_0([0,\infty)\to \R):= \{T:[0,\infty)\to \R;\  T \text{ is c\`adl\`ag}, T(0)=0\}, 
\end{align*}
 induced by the inverse Gaussian subordinator (e.g. \cite{A 09})    
\begin{align*}
U^m(t):=  \inf\{s>0; B^1(s)+ms=t\}, \qquad t \geq 0.  
\end{align*}  
Namely, for the Borel set $E$ in $D_0^{(1)}$, $\nu^m$ is defined by   
$\nu^m(E):=\mu^1 \big(  U^m(\cdot)   \in E\big)$. 
 Here  $B^1$ is $1$-dimensional standard Brownian motion
and  $\mu^1$ is  the $1$-dimensional Wiener measure.       
 Then   $T$  is a    subordinator 
with respect to $\nu^m$. $S_{3,A,V}(x,t;B,T)$ is a complex-valued semimartingale given           by 
\begin{align}\label{semi 3}
S_{3,A,V}    (x,t;B,T)&=i\bigg[ \int_0^{T(t)} A(x+B(s))\cdot dB(s)+\frac{1}{2}\int_0^{T(t)} (\nabla \cdot A) (x+B(s))ds \bigg]\nonumber \\
&\qquad  +\int_0^tV(x+B(T(s)))ds. 
\end{align}

Now, for the probability distributions of $X(t)$ and $B\big(T(t)\big)$,   note that 
\begin{align}\label{probability distribution}
\lambda^m(X(t)\in dy)=\Big(\mu \times \nu^m\Big)\Big(B\big(T(t)\big)  \in dy\Big)=k_0^m(y, t) dy. 
\end{align}
Here $k_0^m(y,t)$ is the integral kernel of the semigroup $e^{-t[\sqrt{-\Delta + m^2}-m]}$ and  has the  explicit expression
\begin{align}\label{eq:eq kernel}
     k_{0}^{m}(y,t) 
 =   \begin{cases}   
\displaystyle 2 \left( \frac{m}{2\pi}\right)^{\frac{d+1}{2}}\frac{te^{mt}  K_{\frac{d+1}{2}}\big(m(|y|^{2}+t^{2})^{\frac{1}{2}}\big   )}{(|y|^{2}+t^{2})^{\frac{d+1}{4}}},  \quad & m>0,
  \\  
 \displaystyle \frac{\Gamma (\tfrac{d+1}{2})}{\pi^{\frac{d+1}{2}}} \frac{t}{(|y|^{2}+t^{2})^{\frac{d+1}{2}}}, & 
m=0. 
 \end{cases}
\end{align}
Expressions \Equref{probability distribution} and \Equref{characteristic f.t. of lambda} imply that
\begin{align}\label{eq;eq fourier transformation of kernel}
\widehat{k_0^m(\cdot, t)}(\xi)=e^{-t[\sqrt{\xi^2+m^2}-m]},\quad \xi \in \R^d. %\ t  > 0.
\end{align}
Here,     
for $\varphi \in \mathcal{S}(\R^{d})$,  we define  the Fourier transform of $\varphi$  by $\widehat{\varphi}(\xi):=\int_{\R^d} e^{-iy \cdot \xi}\varphi(y)dy$.

\begin{remark} Under  the assumption in  Theorem 1.3 (i) (resp. (ii)), 
$H_{1,A}^m-m+V$, $H_{2,A}^m-m+V$ (resp. $H_{3,A}^m-m+V$) can be realized as  nonnegative selfadjoint operators in $L^2(\R^d)$ through the  quadratic  forms (\cite{I and T 93}, \cite{I 2013}). Then each term in  $S^m_{1,A,V}(x,t   ;X)$, $S^m_{2,A,V}(x,t;X)$ (resp. $S_{3,A,V}(x,t;B,T)$) is well-defined $\lambda^m$-a.s. (resp.  $\mu \times \nu^m$-a.s.)     
\end{remark}

\begin{remark} 
Under  the assumption in  Theorem 1.2 (i) (resp. (ii)), the maps 
$x \mapsto S_{1,A,V}^m(x,t;X)$, $x \mapsto S_{2,A,V}^m(x,t;X)$   (resp. $x \mapsto S_{3,A,V}(x,t; B, T)$)    are continuous $\lambda^m$-a.s. (resp.  $\mu \times \nu^m$-a.s.) Especially, then   the third  terms  in  $S_{1,A,V}^m(x,t;X)$, $S_{2,A,V}^m(x,t;X)$  (principal value integrals) are equal  to 
\begin{align*}
& \int_{0}^{t} ds   \int_{0<|y|<1} \left[ A(x+X(s)+\tfrac{1}{2}y) -A(x+X(s))\right]   \cdot y n^{m}(dy), \\ 
&  \int_{0}^{t} ds \int_{0<|y|<1} \left(   \int_0^1 A(x+X(s)+\theta y)  d\theta-A(x+X(s))  \right) \cdot y n^{m}(dy),
\end{align*}
respectively, 
since  $\int_{0<|y|<1}|y|^{1+\alpha}n^m(dy)< \infty$ and $n^m(y)$ is rotatinally invariant. %We do not know whether the assumptions for $A$, $V$ in Theorem 1.2 are the weakest condition for it  to  hold that $x \mapsto S_{1,A,V}^m(x,t;X)$, $x \mapsto S_{2,A,V}^m(x,t;X)$,   $x \mapsto S_{3,A,V} (x,t;   B, T)$ are continuous. 
 \end{remark}

\begin{remark} For the density of L\'evy measure $n^m(dy)$ and the integral kernel $k_0^m(y,t)$  of the semigroup   $e^{-t[\sqrt{-\Delta+m^2}-m]}$, it holds that 
\begin{align*}
n^m(y) \uparrow n^0(y),\  k_0^m(y,t)\uparrow k_0^0(y,t)\quad \text{as }\ m \downarrow 0.
\end{align*}
In fact (e.g. \cite[(21), p.79]{E 53}), since $\frac{d}{d\tau}  \Big(\tau^{\frac{d+1}{2}}K_{\frac{d+1}{2}}
 (\tau)\Big)= -  \tau^{\frac{d+1}{2}}K_{\frac{d-1}{2}}(\tau)<0$ and   $\frac{d}{d\tau}       \Big(e^{\tau}\tau^{\frac{d+1}{2}}K_{\frac{d+1}{2}}(\tau)\Big) = e^{\tau}\tau^{\frac{d+1}{2}}\Big(K_{\frac{d+1}{2}}(\tau)-K_{\frac{d-1}{2}}(\tau)\Big)<0$ for $\tau>0$,   the functions $\tau \mapsto \tau^{\frac{d+1}{2}}K_{\frac{d+1}{2}}(\tau)$ and $\tau \mapsto e^{\tau}\tau^{\frac{d+1}{2}}K_{\frac{d+1}{2}}(\tau)$ are strictly decreasing. Therefore we have %(\cite[(41), (42), (43), p.10]{E 53})
\begin{align*}
\tau^{\frac{d+1}{2}}K_{\frac{d+1}{2}}(\tau),\  e^{\tau}\tau^{\frac{d+1}{2}}K_{\frac{d+1}{2}}(\tau) \uparrow 2^{\frac{d-1}{2}}\Gamma(\tfrac{d+1}{2})\quad  \text{as }\  \tau \downarrow 0.\end{align*}  Then  it follows from 
\Equref{eq:eq levy measure}, \Equref{eq:eq kernel} that as $m \downarrow 0$, 
\begin{align*}
n^{m}(y)&=     
  2\left (\frac{1}{2\pi}\right)^{\frac{d+1}{2}}  \frac{ ( m|y|)^{\frac{d+1}{2}} K_{\frac{d+1}{2}}(m|y|)}{|y|^{d+1}}\uparrow n^0(y),\\ 
 k_0^m(y,t)& = 2 \left( \frac{1}{2\pi}\right)^{\frac{d+1}{2}}     
t e^{-m[(y^2+t^2)^{1/2}-t]}      \\
& \quad \times \frac{  e^{m(y^2+t^2)^{1/2}}     \Big( m(|y|^{2}+t^{2})^{\frac{1}{2}}\Big)^{\frac{d+1}{2}}  K_{\frac{d+1}{2}}\big(m(|y|^{2}+t^{2})^{\frac{1}{2}}\big   )}{(|y|^{2}+t^{2})^{\frac{d+1}{2}}}\uparrow   k^0_0(y,t). 
\end{align*}
\end{remark}

\section{Proof of Theorem 1.1}% and Theorem 3.2.}  

The proof of  claim  (i) of Theorem 1.1 is given in  \cite{I and M 2014}. So  we  prove   only  claim  (ii).  
To this end,  we have  to verify   the following  three facts (\cite[Theorem 13.5]{B 99}):\\
(a) The finite dimensional distributions with respect to $\nu^{m}$  weakly   converge   to those with respect to $\nu^{0}$ as $m \downarrow 0$.\\  
(b) For each  $t > 0$,  the probability measure 
$\nu^0( T(t)-T(s)\in dy )$  weakly converges                                                                       to  Dirac measure concentrated at the point $0\in \R$ as $s \uparrow t$. \\ 
(c) There exist  $\alpha > 0$ and $\beta >1 $, and a nondecreasing continuous function $F$ on $[0, \infty)$  such that for $m>  0$,  $0\leq r \leq s \leq t <\infty$,  $a>0$,  
\begin{align*}
 & \nu^m\Big(\big(T(s)-T(r)\big)\wedge  \big(T(t)-T(s)\big) \geq a \Big) \leq \frac{1}{a^{\alpha}}\Big[F(t)-F(r)\Big]^{\beta}.%\\ & \qquad \qquad \qquad \qquad \quad \qquad \qquad \qquad  m>  0,\  0\leq r \leq s \leq t <\infty, \ a>0. 
\end{align*}
{\it Proof.} To prove (a),    we note that the L\'evy exponent  of $\nu^m$ is given as follows  (\cite[(4.26)]{I 2013}):  %characteristic function of $T(t)$ with respect to $\nu^m$ is (\cite[(4.26)]{I 2013})
\begin{align*}
  \zeta_m(p)& := \frac{  2\sqrt{2}p^2     }{ [(m^2+\sqrt{m^4+4p^2})^{1/2}+\sqrt{2}m] (m^2+ \sqrt{m^4+4p^2})}\\%\nonumber \\
 &  \quad - \frac{\sqrt{2}p}{(m^2+\sqrt{m^4+4p^2})^{1/2}}i,\quad p\in    \R.
\end{align*}
Here, for $m=0$, $p=0$, we understand 
$\zeta_0(0):= 0$. It is easy to see that
$\zeta_m(p)\to \zeta_0(p)$
as $m \downarrow 0$ for any $p \in \R$. 
Then     (a) follows from this convergence and  independent and stationary increments property of subordinator $T$.

Next,   (b)  follows from the stochastic continuity of subordinator $T$. 

Finally, we prove  (c). Since 
\begin{align}\label{density of subordinator}    
\nu^m(T(t)\in dr)=\frac{t}{\sqrt{2\pi}}\ e^{mt}  r^{-\frac{3}{2}}\exp\left\{-\tfrac{1}{2}(\tfrac{t^2}{r}+m^2r)\right\}dr, \quad r>0,  
\end{align}
 we have for $t> 0$, $a>0$   
\begin{align*}
\nu^m(T(t)\geq a)
& =\frac{t}{\sqrt{2\pi a}}\ e^{mt}\int_1^\infty s^{-\frac{3}{2}}\exp\left\{-\tfrac{1}{2}(\tfrac{t^2}{as}+m^2as)\right\}ds\\
& \leq \frac{t}{\sqrt{2\pi a}}\ e^{mt}  \int_1^\infty s^{-\frac{3}{2} }\exp\left\{-\tfrac{1}{2}\cdot 2\sqrt{\tfrac{t^2}{as}\cdot m^2as}\right\}ds\\
& =  \sqrt{\frac{2}{\pi a }} \ t. 
\end{align*}
From the   independent increments property of subordinator $T$ and the above estimate, we have 
\begin{align*}
  \nu^m\Big(  \big(T(s)-T(r)\big)\wedge  \big(T(t)-T(s)\big)% \\
\geq a \Big)
& \leq  \frac{1}{a} \left(\frac{t}{\sqrt{\pi}}-\frac{r}{\sqrt{\pi}} \right)^2,\quad 0\leq r\leq s\leq t<\infty. 
\end{align*}
Therefore (c) holds for $\alpha:=1$, $\beta:=2$ and $F(t):=\frac{t}{\sqrt{\pi}}.$ 
\hfill$\square$

\section{Change of probability measures}
In this section, we give preliminaries to prove             Theorem 1.2 and Theorem 1.3. %As mentioned in Section 1, 
The main idea is to  change   probability measures on the  right-hand sides of (1.5) and (1.6)        from $\lambda^m$  and $\mu \times \nu^m$ to $\lambda^0$ and $\mu \times \nu^0$, respectively.    More  precisely,         we find   path transformations  $\Phi_m: D_0\ni X \mapsto \Phi_m(X)\in D_0$ and  $\Psi_m:D_0^{(1)}\ni T\mapsto \Psi_m(T) \in D_0^{(1)}$ such that by $\lambda^m=\lambda^0 \Phi_m^{-1}$  and    $\nu^m=\nu^0 \Psi_m^{-1}$,  the right-hand sides of (1.5) and  (1.6) are rewritten as   
\begin{align}
(e^{-t[H_{j, A}^m-m+V]}g)(x)& = E^{\lambda^0}\bigg[e^{-S_{j,A,V}^m(x,t; \Phi_m(X))}g\big(x+\Phi_m(X)(t)\big)\bigg] \quad j=1,2,\\
(e^{-t[H_{3, A}^m-m+V]}g)(x)& = E^{\mu \times \nu^0}\bigg[e^{-S_{3,A,V}(x,t; B, \Psi_m(T))} g\Big(x+B\big(\Psi_m(T)(t)\big)\Big)\bigg].
\end{align}
 Here $\Phi_m(X)(t)$ and $\Psi_m(T)(t)$ are the values of $\Phi_m(X)$ and $\Psi_m(T)$ at time $t$, respectively.  

In fact, in \cite{I and M 2014}, the L\'evy process $\Phi_m(X)$ with respect to $\lambda^0$ has  already been obtained  through   a mapping  $\phi_m:\R^d\setminus \{0\}\to \R^d\setminus \{0\}$ satisfying  
$n^m(dy)=n^0\phi_m^{-1}(dy)$. Namely we have   
  a strictly increasing function $\ell_m:(0,\infty)\to (0,\infty)$ such that 
\begin{align}
 \Phi_m(X)(t) & :=  \int_{0}^{t} \int_{|z| \geq 1}\phi_{m}(z) N_{X}(ds dz)+ \int_{0}^{t}\int_{0< |z| < 1}  \phi_{m}(z) \widetilde{N_{X}^{0}}( ds  dz),\quad \lambda^0{\text-a.s.}, \\
 \phi_m(z)&:=\ell_m^{-1} (|z|)\frac{z}{|z|},\qquad \ell_m(r):=  \frac{2^{\frac{d-1}{2}} \Gamma (\frac{d+1}{2})}{ m^{\frac{d+1}{2}}   \int_{r}^{\infty}    u^{\frac{d-3}{2}} K_{\frac{d+1}{2}}(mu)du }. 
\end{align}
Here defining   $\ell_0(r):=r$ for $m = 0$, we have  $\phi_0(z)=z,$ $\Phi_0(X)=X$. Since $\ell_m(r) \downarrow r$ as $m\downarrow 0$ (\cite[Proposition 1 (ii)]{I and M 2014}),   we have
\begin{align}
\phi_m(z)\to z, \; |\phi_m(z)|=\ell_m^{-1}(|z|) \uparrow |z|\quad \text{as }\ m \downarrow 0.  
\end{align}
Therefore we obtain the following proposition (\cite[Proposition 2]{I and M 2014}): 
\begin{proposition}
For every sequence $\{m\}$ with $m\downarrow 0$, there exists a subsequence $\{m'\}$ such that
\begin{align*}
\sup_{t \leq t_0 }  | \Phi_{m'}(X)(t)-X(t)|\to 0 \quad \text{as }  \  m'   \downarrow 0, \ \lambda^0\text{-a.s.}
\end{align*}
\end{proposition}

Next, by an analogous argument used to obtain $\Phi_m$ in \cite{I and M 2014}, we will    find $\Psi_m:D_0^{(1)}\to D_0^{(1)}$ such that $\nu^m=\nu^0\Psi_m^{-1}$. Let $\sigma^m(dr)$ be the L\'evy measure of  subordinator $T$ with respect to $\nu^m$.  It is known that  
\begin{align*}
\int_0^\infty \frac{f(r)}{t}\nu^m(T(t)\in dr) \to \int_0^\infty f(r)  \sigma^m(dr) \quad  \text{as }\ t \downarrow 0,  
\end{align*}
for any bounded continuous function $f\colon (0,\infty)\to \R$ vanishing in a neighborhood of the origin  (\cite[(6.4.11)]{N 2000}). 
% $\frac{\nu^m(T(t)\in dr)}{t} \to o^m(dr)$ as $t \downarrow 0$, 
Then we have by \Equref{density of  subordinator}  that 
\begin{align*}
\sigma^m(dr)&= \frac{1}{\sqrt{2\pi} }   r^{-\frac{3}{2}}e^{-\frac{1}{2}m^2r}
dr =:\sigma^m(r)dr. 
\end{align*}
Now,  we will determine $\psi_m \colon (0,\infty) \to (0, \infty)$ in such a way that (i) $\sigma^m(dr)=\sigma^0\psi_m^{-1}(dr)$, (ii) $\psi_m\in C^1((0,\infty);(0,\infty))$, (iii) $\psi_m$ is bijective and  (iv) $\psi_m'(r) \ne 0$ for 
all $r>0$. For any  Borel set $U$  in $(0,\infty)$,   we have  
\begin{align*}
\sigma^m(U)=\int_U \sigma^m(r)dr,\quad 
\sigma^0 \psi_m^{-1}(U)= \int_U \sigma^0(\psi_m^{-1}(r)) (\psi_m^{-1})'(r) dr. 
\end{align*}
Therefore we have  
$\sigma^m(r)= \sigma^0(\psi_m^{-1}(r)) (\psi_m^{-1})'(r)$  a.s. $r>0$, 
and hence 
\begin{align*}
r^{-\frac{3}{2}}e^{-\frac{1}{2}m^2r}= (\psi_m^{-1}(r))^{-\frac{3}{2}}(\psi_m^{-1})'(r)       \quad \text{a.s. }r>0.
\end{align*}
We solve this differential equation under boundary condition $\psi_m^{-1}(\infty)=\infty$   to get 
\begin{align*}
\psi_m^{-1}(r)= \frac{4}{   \left(  \int_r^\infty u^{-\frac{3}{2}}e^{ -\frac{1}{2}m^2u }du    \right)^2}.
\end{align*}
Since $\psi_m^{-1}(r)\downarrow r$,  we have 
$\psi_m(r)\uparrow r$     as $m \downarrow 0$. 
For $m=0$, we put $\psi_0(r):=r$.  Thus we  determined $\psi_m$. Next,  by noting  $\int_0^\infty r \ \sigma^0(dr)< \infty$,  we define subordinator $\Psi_m(T)$  % for $m>0$ 
with respect to $\nu^0$ by
 \begin{align}\label{path of transformed subordinator}
\Psi_m(T)(t)=\int_0^t\int_0^\infty \psi_m(r)  N_T(dsdr)\quad \nu^0\text{-a.s.} 
\end{align} for $m \geq 0$. Here  
\begin{align*}
N_T(dsdr):=\# \{u\in ds; T(u)-T(u-)\in dr\},\quad s >0, r>0.
\end{align*} 
It is trivial that $\Psi_0(T)=T$,  $0\leq\Psi_m(T)(t) \leq T(t)$. It can been seen that  $\nu^m=\nu^0\Psi_m^{-1}$. In fact, for $p_1, \ldots,  p_k\in\R$,  $0=s_0<s_1<\cdots<  s_k<\infty$, $k\in \mathbf{N}$, we have  by  \cite[Theorem 1.3.15]{A 09} and the relation $\sigma^m(dr)=\sigma^0\psi_m^{-1}(dr)$ that   
\begin{align*}
E^{\nu^m}\left[e^{i      \sum_{j=1}^k p_{j} T(s_{j})}  \right] 
&= \prod_{j=1}^{k}          \exp\left\{ (s_j- s_{j-1}) \int_0^\infty (e^{i(p_j+\cdots+ p_k)r}-1)\sigma^m(dr)    \right\}\\
&= \prod_{j=1}^{k}          \exp\left\{ (s_j- s_{j-1}) \int_0^\infty (e^{i(p_j+\cdots+  p_k)\psi_m(r) }-1)\sigma^0(dr)    \right\}\\
&= E^{\nu^0}\left[ e^{i    \sum_{j=1}^k p_{j} \Psi_m(T)(s_{j})}  \right].  
\end{align*}

Now,  
we can get also the following proposition for $\Psi_m(T)$ corresponding to Proposition 4.1 for $\Phi_m(X)$. Its  proof  is easy since 
$\psi_m(r)\uparrow r$     as $m \downarrow 0$:
\begin{proposition}
For every sequence $\{m\}$ with $m \downarrow 0$, we have 
\begin{align*}
(0\leq )\ \sup_{t \leq t_0} (T(t)-\Psi_m(T)(t))\to 0\quad \text{as }\ m\downarrow 0, \nu^0\text{-a.s.}
\end{align*}
\end{proposition}

\vspace{2mm}

It is to be noted that we need to take a subsequence in Proposition 4.1 but not in Proposition 4.2.

\section{Proof of Theorem 1.2}
In this section, we prove Theorem 1.2. 
First we prove   two  key  lemmas. 
\begin{lemma} 
Assume that $0\leq V\in C(\R^d;\R)$ and         $0<t_0<\infty$, $0<R      < \infty$. \\     
{\rm(i)} If $A$ is locally $\alpha$-H\"older continuous  $(0<\alpha \leq 1)$, then it holds that for $j=1,2$
\begin{align*}
E^{\lambda^0}\left[
\left| e^{-S_{j,A,V}^{m}(x , t;  \Phi_{m}(X))}  - e^{-S^0_{j,A,V}( x,t; X)} \right| 
\right]
\to 0 \quad \text{as }\     m \downarrow 0,   \end{align*} 
uniformly on   $t \leq t_0, |x|<R.$\\
{\rm(ii)}  If $A \in C^1(\R^d;\R^d)$, then it holds that      
\begin{align*}     
 E^{\mu\times \nu^0}\left[%\mathbf{1}_{  t<\sigma_k (\Phi_m(X))\wedge  \sigma_k(X)}                                             
\bigg| e^{-S_{3,A,V}(x , t; B,  \Psi_{m}(T))}  - e^{-S_{3,A,V}  ( x,t;B,  T)} \bigg| 
\right]
\to 0\quad 
\text{as }\  
m \downarrow 0,    \end{align*} 
uniformly on   $t \leq t_0, |x|<R.$
\end{lemma}
\noindent 
{\it Proof.} (i) First we prove claim  (i) for $j=1$.  
By  $N_{\Phi_m(X)}(dsdy)=N_{X}(ds\phi_m^{-1}(dy))$,         
it follows from  \Equref{semi 1} that  
\begin{align}\label{semimartingale 1}
S^{m}_{1,A,V}       (x,t; \Phi_m(X))&  =   i  \bigg[  \int_{0}^{t}\int_{|z| \geq 1} A(x+\Phi_m(X)(s-)+\tfrac{1}{2}\phi_m(z)   ) \cdot \phi_m(z) \ N_{X}(dsdz) \nonumber \\
 & \quad \quad + \int_{0}^{t}  \int_{0<|z| < 1}   A(x+\Phi_m(X)(s-)+\tfrac{1}{2}\phi_m(z)) \cdot \phi_m(z) \      \widetilde{N_{X}^{0}}(dsdz)\nonumber \\
 & \quad\quad + \int_{0}^{t} ds \text{ p.v.} \int_{0<|z|<1}  A(x+\Phi_m(X)(s)+\tfrac{1}{2}\phi_m(z)) \cdot \phi_m(z) \  n^{0}(dz)  \bigg]  \nonumber   \\
& \quad \quad + \int_{0}^{t} V(x+\Phi_m(X) (s))\ ds\nonumber \\
&=: i\Big[S_{1,A}^m(x,t;X)+S_{2,A}^m(x,t;X)+S_{3,A}^m(x,t;X)\Big]   +S_{4,V}^m(x,t;X).
\end{align}
Then we have 
\begin{align} \label{eq:eq main estimate}
 & \sup_{t \leq t_0,  |x|<R } E^{\lambda^0}\left[
           \bigg| e^{-S^{m}_{1,A,V}(x,t; \Phi_{m}(X))}  - e^{-S^{0}_{1,A,V} (x,t;  X)} \right| 
\bigg]\nonumber \\
& \leq E^{\lambda^0}\bigg[\sup_{t \leq t_0, |x|<R  }
 \Big|e^{-iS_{1,A}^{m}( x, t; X)}- e^{-i S_{1,A}^0(x, t; X)}  \Big|
\bigg]       +  \sup_{|x|<R }E^{\lambda^0}\bigg[ \sup_{t \leq t_0} 
\Big|e^{-iS_{2,A}^{m}(x, t; X)}- e^{-i S_{2,A}^0( x, t; X)}  \Big|\bigg] \nonumber \\
    &  \quad +E^{\lambda^0}\bigg[  \sup_{t \leq t_0, |x|<R } 
\Big|e^{-iS_{3,A}^{m}( x, t; X)}- e^{-i S_{3,A}^0(x, t; X)}  \Big|
 \bigg] +E^{\lambda^0}\bigg[\sup_{t \leq t_0,|x|<R } \Big|e^{-S_{4,V}^{m}(x, t; X)}- e^{-S_{4,V }^0( x, t; X)}\Big|
 \bigg]\nonumber \\
&=: E^{\lambda^0}\big[I_1^m(X)\big]   +   \sup_{|x|<R} E^{\lambda^0}\big[I_2^m(x; X)\big]    +E^{\lambda^0}\big[I_3^m(X)\big]  +E^{\lambda^0}\big[I_4^m(X)\big].  
\end{align}  
We   now  show that each term in  the last member    of 
\Equref{eq:eq main estimate}       converges  to zero as $m \downarrow  0$. To this end,  we note that  $I_1^m(X)$, $I_3^m(X)$ and $I_4^m(X)$ are less than or equal to   $2$.  Let $\{m\}$ be a sequence with $m\downarrow 0$ and $\{m'\}$ any subsequence of $\{m\}$. Then,       by  Proposition 4.1, there exists a subsequence $\{m''\}$ of $\{m'\}$ such that $\displaystyle \sup_{t \leq t_{0}}|\Phi_{m''} (X)-X(t)|\to 0$ as $m''\downarrow 0$, $\lambda^0$-a.s.

{\it For the first term       of  \Equref{eq:eq main estimate}:}      By the definition of $N_{X}(dsdy)$, we have    
\begin{align*}
  &  S^{m''}_{1,A}(x,t; X) -S^{0}_{1,A}(x,t; X)  \\
&=  \sum_{s\leq t}  \mathbf{1}_{|X(s)-X(s-)|\geq 1}\bigg[ \big( A(x+\Phi_{m''}(X)(s-) + \tfrac{1}{2}\phi_{m''}(X(s)-X(s-))) \\
&\qquad \qquad \qquad    \qquad -  A(x+X(s-)+\tfrac{1}{2}(X(s)-X(s-)))\big)\cdot  \phi_{m''} (X(s)-X(s-)) \\ 
&  \qquad \qquad +A(x+\tfrac{1}{2}(X(s)+ X(s-)))\cdot   \big(\phi_{m''}(X(s)-X(s-))- (X(s)-X(s-))\big)  \bigg],
\end{align*}
which   is a              finite  sum  (e.g. \cite[p.122]{B 99}). 
Then    we have     
\begin{align*}
 I_1^{m''}(X)\leq 
& \sum_{s\leq t_0}  \mathbf{1}_{|X(s)-X(s-)|\geq 1}\bigg[       C_1(X ) \Big(|\Phi_{m''}  (X)(s-)-X(s-)|\\
&\qquad \qquad \qquad \qquad \qquad +\tfrac{1}{2}|\phi_{m''}(X(s)-X(s-))-(X(s)-X(s-))|\Big)^\alpha C_2(X)       \\
&\quad   +  C_3(X) 
\big|\phi_{m''}(X(s)-X(s-))-(X(s)-X(s-))\big|    \bigg], 
\end{align*}
since  $A$ is locally  $\alpha$-H\"older continuous and so locally  bounded.  
Here  $C_1(X), C_2(X), C_3(X)$ are  constants  depending on $X$. Since $\phi_{m^{''}} (z)\to z$,  the above  sum  converges to zero as $m'' \downarrow 0$, $\lambda^0$-a.s.     
   Hence $\displaystyle E^{\lambda^0}[I_1^{m''}(X)]$ converges to  zero as $m'' \downarrow 0$.

{\it For the          second  term of  \Equref{eq:eq main estimate}:}                                                                                                                                       %Finally, as for the second term, 
First, for $k\in \mathbf{N}$, let  $\sigma_k(X)$ be  the  hitting time %for $\{y;|y|>k\}$,  namely  % the $\{\mathcal{F}(t)\}_{t \geq 0}$-stopping time 
defined by 
\begin{align}
\sigma_k(X):  =   
\inf\{s>0; |X(s-)|> k\}. 
\end{align}
Here we understand $\inf \emptyset :=\infty$  if the  set $\{s>0; |X(s-)|> k\}$ is empty. 
Then  it holds that 
$\sigma_{k}(X) \to  \infty$ as $k \to \infty$ and  
$|X(s-)| \leq k$ for   $0<s\leq \sigma_k(X)$.   
From  the relation $\displaystyle 
\int_0^t=\int_0^{t \wedge \sigma_k(X) \wedge \sigma_k(\Phi_{m''}(X))} + \int_{t \wedge \sigma_k(X) \wedge \sigma_k(\Phi_{m''}(X))}^t$ 
and Doob's martingale inequality, we have %$I_2^m(x;X)$ is less than or equal to  
\begin{align}\label{martingale term} 
E^{\lambda^0}[I_2^{m''}(x;X)]    & \leq        2E^{\lambda^0} \bigg[\int_0^{t_0      \wedge \sigma_k(X)\wedge\sigma_k(\Phi_{m''}(X))}ds\nonumber \\
& \qquad \qquad \times \int_{0<|z|<1}\big|
A(x+\Phi_{m''} (X)(s-)+\tfrac{1}{2}\phi_{m''} (z)) \cdot \phi_{m''}(z)\nonumber \\
&\qquad \qquad  \qquad \qquad - A(x+X(s-)+\tfrac{1}{2}z) \cdot     z\big|^2  n^0(dz)    \bigg]^{1/2}\nonumber \\
&  \quad + 2\lambda^0(\sigma_k(X)<  t_0)+2\lambda^0(\sigma_k(\Phi_{m''}(X))<t_0). 
\end{align}
Note that  
\begin{align*}
&A(x+\Phi_{m''}(X)(s-)+\tfrac{1}{2}\phi_{m''}(z)) \cdot \phi_{m''}(z)  - A(x+X(s-)+\tfrac{1}{2}z) \cdot     z\\
& =\Big[A(x+\Phi_{m''}(X)(s-)+\tfrac{1}{2}\phi_{m''}(z))  - A(x+X(s-)+\tfrac{1}{2}z)\Big] \cdot \phi_{m''}(z)\\
& \quad + A(x+X(s-)+\tfrac{1}{2}z) \cdot (\phi_{m''} (z)-z).
\end{align*}  
Since $(a+b)^2\leq 2a^2+2b^2$ ($a,b\in\R$) and  
$|\phi_{m''}(z)|\leq |z|$, 
  the first term on the right-hand side of \Equref{martingale term} is less  than or equal to 
\begin{align*}
& 2 \sqrt{2} \Bigg(E^{\lambda^0}\Bigg[\int_0^{t_0} ds\int_{0<|z|<1} \sup_{\substack{|w|,|w'| < R+k+\frac{1}{2};\\ |w-w'|\leq |\Phi_{m''}(X)(s-)-X(s-) | + \frac{1}{2}|\phi_{m''}(z)-z| } } |A(w)-A(w')|^2 |z|^2n^0(dz)\Bigg]\\
& \quad\quad +     t_0 \sup_{|w|<R+k+\frac{1}{2}}|A(w)|^2 \int_{0<|z|<1}|\phi_{m''} (z)-z|^2n^0(dz)\Bigg)^{1/2}.
\end{align*}
This  converges to zero as $m'' \downarrow 0$ since  $A$ is locally uniformly continuous.  Then we have by \Equref{martingale term} that  
\begin{align*}
 \displaystyle \limsup_{m''\downarrow 0} \sup_{|x|<R} E^{\lambda^0}\bigg[ I_2^{m''}(x;X)\bigg] & \leq 2\lambda^0(\sigma_k(X)<  t_0)+2\lambda^0\left(\limsup_{m''\downarrow 0}\Big\{ \sigma_k(\Phi_{m''}(X))<t_0\Big\}\right)\\ 
&  \leq  2\lambda^0(\sigma_k(X)<  t_0)+2 \lambda^0(\sigma_{k-  1}(X) < t_0),
\end{align*} 
which converges to zero as $k\to \infty$ because  of  $\sigma_k(X)\to \infty$.  Hence $\displaystyle \sup_{|x|<R} E^{\lambda^0}\bigg[ I_2^{m''}(x;X)\bigg]$ converges to zero  as $m'' \downarrow 0$.

{\it For  the third   term of  \Equref{eq:eq main estimate}:  }                                                                                                                                      
In view of Remark 2.2, we have 
\begin{align*}
& S^{m''}_{3,A}( x,t; X)-S^{0}_{3,A}( x,t; X)  \\
&=  \int_0^t ds \int_{0<|z|<1}      \left[A(x+\Phi_{m''} (X) (s)+\tfrac{1}{2}\phi_{m''}(z)  )-A(x+\Phi_{m''}(X)(s))\right] \cdot (\phi_{m''} (z)-z)  n^{0}(dz)\\
& \quad + \int_0^t ds \int_{0<|z|<1}   \bigg[   \Big[A(x+\Phi_{m''} (X) (s)+\tfrac{1}{2}\phi_{m''}(z)  )  -A(x+\Phi_{m''}(X)(s))\Big] \\
&\quad \quad \qquad \qquad\qquad  -\Big[A(x+X(s)+\tfrac{1}{2}z)-A(x+X(s))\Big]\bigg]\cdot z  n^{0}(dz).
\end{align*}
It follows from the above expression and  the  local  $\alpha$-H\"older continuity of $A$ that 
\begin{align*}
 I_3^{m''}(X)&      \leq     t_0 C(X )  \int_{0<|z|<1} \left(\frac{1}{2}|\phi_{m''} (z)|\right)^\alpha   
|\phi_{m''} (z)-z|  n^{0}(dz)\nonumber \\
& \quad + \int_0^{t_0} ds \int_{0<|z|<1}  \sup_{|x|<R}  \bigg|   \Big[A(x+\Phi_{m''} (X) (s)+\tfrac{1}{2}\phi_{m''}(z)  )  -A(x+\Phi_{m''}(X)(s))\Big] \nonumber \\
&\quad \quad \qquad \qquad\qquad  \qquad \qquad  -\Big[A(x+X(s)+\tfrac{1}{2}z)-A(x+X(s))\Big]\bigg||z
|  n^{0}(dz)\\
&=: t_0 C(X) \int_{0<|z|<1} J_1^{m''}(z) n^0(dz)+ \int_0^{t_0}ds \int_{0<|z|<1} J_{2}^{m''}(s,z;X)n^0(dz). 
\end{align*}
Here $C(X)$ is a constant depending  on $X$.  For $J_1^{m''}(z)$, since $|\phi_{m''}(z)|\leq |z|$ and \\ $\int_{0<|z|<1}     |z|^{1+\alpha} n^0(dz)<\infty$,  $\displaystyle \int_{0<|z|<1} J_1^{m''}(z) n^0(dz)$    converges to zero as $m'' \downarrow 0$.  For $J_{2}^{m''}(s,z;X)$,  it is easy to see that $J_2^{m''}(s,z;X)$ converges to zero as $m'' \downarrow 0$ $\lambda^0$-a.s. for fixed $s$ and $z$.      
On the other hand,    we have    
\begin{align*}
 J_{2}^{m''}(s,z;X)\leq C(X ) \left( (\tfrac{1}{2}|\phi_{m''}(z)|)^\alpha+(\tfrac{1}{2}|z|)^\alpha \right) |z|\leq 
 C(X )\tfrac{1}{2^{\alpha-1}}|z|^{1+\alpha}.
\end{align*} 
Therefore  $\displaystyle \int_0^{t_0} ds \int_{0<|z|<1} J_{2}^{m''}(s,z)n^0(dz)$ converges to zero as   $m''\downarrow 0$ $\lambda^0$-a.s.  Hence $E^{\lambda^0}[I_3^{m''}  (X)]$ converges to zero  as $m'' \downarrow 0$.

{\it For  the fourth term of              \Equref{eq:eq main estimate}:}   Note that  $V \in C(\R^d;\R)$     is locally  uniformly continuous. Then we have  
\begin{align*}
 I^{m''}_4(X)\leq 
 \int_{0}^{t_0 } \sup_{|x|<R} \Big|V\big(x+ \Phi_{m''}(X)(s)\big)- V\big(x+ X(s)\big)\Big|
ds,\end{align*} which converges  to zero as $m'' \downarrow 0$ $\lambda^0$-a.s.. Hence  $E^{\lambda^0}[ I^{m''}_4(X)]$  converges to zero  as $m'' \downarrow 0$.

Thus we have seen that the four terms of (5.2) converges to zero as $m \downarrow 0$, which shows claim (i) for $j=1$.

The convergence  for $j=2$  can be proved in the same way as for $j=1$ above.  In fact, we have only to replace $A(x+\Phi_{m''}(X)(s-)+\frac{1}{2}\phi_{m''}(z))\cdot \phi_{m''}(z)$ and $A(x+X(s-)+\frac{1}{2}z)\cdot z$  by $\left(\int _0^1A(x+\Phi_{m''}(X)(s-)+\theta \phi_{m''}(z))d\theta \right)\cdot \phi_{m''} (z)$ and $\left(\int _0^1A(x+X(s-)+\theta z)d\theta \right)\cdot z$,  respectively.  
This shows claim (i) for $j=2$, ending the proof of claim (i).

\vspace{5mm}

\noindent
(ii) By \Equref{semi 3}, we obtain 
\begin{align}\label{semimartingale 3}
S_{3,A,V} (x,  t; B, \Psi_m(T))
& = i\bigg[\int_0^{\Psi_m(T)(t)}A(x+B(s))\cdot dB(s) +\frac{1}{2} \int_0^{\Psi_m(T)(t)} (\nabla\cdot A)( x+B(s)) ds\bigg] \nonumber   \\
&\quad \quad +  \int_0^t V\Big(x+B\big(\Psi_m(T)(s)\big)\Big) ds\nonumber \\
&=: i\Big[ S_{1,A}^m(x,t; B,T)+S_{2,A}^m(x,t; B,T)\Big]+S_{3,V}^m(x,t;B,T). 
\end{align}   
Then      we   have
\begin{align}\label{main estimate 2}
 &\sup_{t \leq t_0,|x|<R } E^{\mu\times \nu^0}\bigg[%\mathbf{1}_{  t<\sigma_k (\Phi_m(X))\wedge  \sigma_k(X)}                                             
\Big| e^{-S_{3,A,V}(x , t; B,  \Psi_{m}(T))}  - e^{-S_{3,A,V}( x,t;B,  T)} \Big| %|g(\cdot +X(t))|  %\left|g(\cdot+X(t))\right|   
\bigg] \nonumber \\%
& \leq E^{\nu^0} \left[\ \sup_{|x|<R}  E^{\mu}\left[\sup_{t\leq t_0}\left| e^{-iS_{1,A}^m(x , t; B,  T)}  - e^{-i S_{1,A}^0( x,t;B,  T)} \right|  \right] \right]\nonumber \\ %     ]
& \quad + E^{\mu\times \nu^0}\left[\sup_{t\leq t_0,|x|<R}\left| e^{-iS_{2,A}^m(x , t; B,  T)}  - e^{-i S_{ 2,A}^0( x,t;B,  T)} \right|  \right]\nonumber \\ %     ]
& \quad + E^{\mu\times \nu^0}\left[\sup_{t\leq t_0,|x|<R}\left| e^{-S_{3,V}^m (x , t; B,  T)}  - e^{-S_{3,V}^0( x,t;B,  T)} \right|  \right]\nonumber \\
& =: E^{\nu^0}\left[\sup_{|x|<R} E^{\mu}\big[I_1^m(x; B,T)\big] \right]+ E^{\mu \times \nu^0}\big[I_2^m(B,T)\big]+ E^{\mu \times \nu^0}\big[I_3^m(B,T)\big]. 
\end{align}We  now    show each  term in the last member   of 
\Equref{main estimate 2}       converges  to zero as $m \downarrow  0$.

{\it For  the first  term                  of              \Equref{main estimate 2}:} 
Note that $\Psi_m(T)(t)\leq T(t)$. From the relations 
\begin{align*}  
&S_{1,A}^m(x , t; B,T)=\int_0^{\Psi_m(T)(t) \wedge \sigma_k(B)}+ \int_{\Psi_m(T)(t) \wedge \sigma_k(B)}^{\Psi_m(T)(t)},\\
&S_{1,A}^0(x , t; B,   T)=\int_0^{T(t)  \wedge \sigma_k(B)}+ \int_{T(t)  \wedge \sigma_k(B)}^{T(t) },
\end{align*}
and  Doob's martingale inequality,  we have 
\begin{align}\label{Brownian part estimate}
E^{\mu}\big[I_1^m(x; B,T)\big]& \leq
 2E^\mu\left[   \int_{\Psi_m(T)(t_0)\wedge  \sigma_k(B)}^{T(t_0) \wedge \sigma_k(B)} | A(x+B(s))|^2ds\right]^{1/2}  \nonumber\\
&\quad + 2\mu(\sigma_k(B) < \Psi_m(T)(t_0))+2\mu(\sigma_k(B)<T(t_0)). 
\end{align}
From  Proposition 4.2 and the fact that $A$ is locally bounded, we have   
\begin{align*}
  \int_{\Psi_m(T)(t_0)\wedge  \sigma_k(B)}^{T(t_0) \wedge \sigma_k(B)}\sup_{|x|<R} | A(x+B(s))|^2ds    
\begin{cases}
&\to 0\quad \text{as }\ m \downarrow 0\  \nu^0\text{-a.s.}, \\
&\leq \displaystyle \sup_{|z|< R+ k} |A(z)|^2 \ T(t_0)<\infty.   
\end{cases}
\end{align*}
By the above   and \Equref{Brownian part estimate}, we have
\begin{align*}
\limsup_{m\downarrow 0}  \sup_{|x|<R}  E^{\mu}\big[I_1^m(x; B,T)\big]
 \leq 
 4\mu(\sigma_k(B)<T(t_0)), % \to 0, %\ \text{as }\ m \downarrow 0,\\
\end{align*}
which converges to zero as $k \to \infty$. 
%which converges to zero as $k\to \infty$. 
Hence  $\displaystyle E^{\nu^0}\bigg[\sup_{|x|<R}  E^{\mu}\big[I_1^m(x; B,T)\big]\bigg]$  converges to zero  as $m \downarrow 0$. 

{\it   For the  second    and  third  terms          of              \Equref{main estimate 2}:}  Note that $\nabla \cdot A$ and $V$ are  locally bounded. 
Then we have 
\begin{align*}
&I_2^m(B,T)\leq  C(B, T)  \ \sup_{t \leq t_0}\big(T(t)-\Psi_m(T)(t)\big),\\ &I_3^m (B,T)\leq   \int_0^{t_0} \sup_{|x|<R} \bigg|V\Big(x+B\big(\Psi_m(T)(s))\Big)-V\Big(x+ B(T(s))\Big) \bigg|ds,
\end{align*}
which    converge to zero as $m\downarrow 0$ $\mu \times \nu^0$-a.s.   Here  $C(B,T)$ is a constant depending on $B$, $T$. 
  Hence  $E^{\mu\times\nu^0}[I_2^m(B,T)]$    and   $E^{\mu\times \nu^0}[I_3^m(B,T)]$                                                                                                                converge to zero  as $m \downarrow 0$.

This shows claim (ii), completing the proof of Lemma 5.1.   \hfill$\square$

\begin{lemma}
$\int_{|y|\geq R}k_0^m(y,t)dy$       converges to zero as $R \to \infty$, uniformly on $m \geq 0$, $t \leq t_0$.   
\end{lemma}
\noindent 
{\it Proof. }    
   Let $\chi$ be a nonnegative $C^{\infty}_{0}$ function  with $0 \leq \chi(y)\leq 1$ in $\R^d$   such that $\chi(y)=1$ if $|y| \leq \frac{1}{2}$ and $\chi(y)=0$ if $|y| \geq 1$. The function $\chi$ satisfies  $\mathbf{1}_{|y|< R} \geq \chi(\frac{y}{R})$ and $\widehat{\chi(\tfrac{\cdot}{R})}(\xi) =R^d   \widehat{\chi}(R\xi)$. Then it follows from Parseval's equality and  \Equref{eq;eq fourier transformation of kernel} that 
\begin{align*}
                                    \int_{|y|\geq R}k_0^m(y,t)dy    &\leq   \int_{\R^d} (1-  \chi(\tfrac{x}{R}))   k_0^m(y,t)dy\\
                                           &= 1- \frac{1}{(2\pi)^d} \int_{\R^d}    \widehat{\chi}(\eta) e^{-t   \big[\sqrt{\frac{\eta^2}{R^2}+m^2}-m\big]  } d\eta\\
& \leq  1- \frac{1}{(2\pi)^d} \int_{\R^d}    \widehat{\chi}(\eta) e^{-t_0   \frac{|\eta|}{R}   } d\eta,
\end{align*}
which converges to zero     as $R \to \infty$. 
This ends  the proof of Lemma 5.2.  \hfill$\square$

\vspace{5mm}

Finally  we prove Theorem 1.2.  First, we show claim (i).  Suppose $g \in C_\infty(\R^d)$ and consider the case   $j=1$ or $2$. Then we have  by (4.1) that 
\begin{align}\label{eq:eq last proof}
\|e^{-t[H^m_{j,A}-m+V]}g- e^{-t[H^0_{j,A}+V]}g\|_\infty &\leq \left\| E^{\lambda^0}\left[\left|e^{-S^m_{j,A,V}(\cdot ,t;\Phi_m(X))}-e^{-S^0_{j,A,V}(\cdot,t;X)} \right| |g(\cdot+X(t))| \right]\right\|_\infty \nonumber \\
& \quad +  E^{\lambda^0}\Big[\| g(\cdot+\Phi_m(X)(t))-g(\cdot+X(t))\|_\infty\Big].
\end{align}
Since  $g$ is uniformly continuous on $\R^d$, the      second term on the right-hand side of \Equref{eq:eq last proof} converges to zero as $m   \downarrow 0$ uniformly on  $t \leq t_0$.   On the other hand,   the first term on the right-hand side of \Equref{eq:eq last proof} is less than or equal to 
\begin{align*}
&   \|g\|_\infty \sup_{|x|<R} E^{\lambda^0}\left[\left|e^{-S^{m}_{j,A,V}(x ,t;\Phi_{m}(X))}-e^{-S^0_{j,A,V}(x,t ;X)} \right|  \right]\vee 2\sup_{|x|\geq R}  E^{\lambda^0}\Big[%|e^{-S_m^{(j)}(x ,t,;\Phi_m(X))}-e^{-S_0^{(j)}(x,t,;X)} | 
\big|g(x+X(t))\big| \Big]
\end{align*}
for $R>0$. Therefore we have  from  Lemma 5.1 (i)  that      
\begin{align*}%\label{eq:eq last proof 2}
&\limsup_{m   \downarrow 0} \sup_{t\leq t_0} \|e^{-t[H^m_{j,A}-m+V]}g- e^{-t[H^0_{j,A}+V]}g\|_\infty \\
& \leq   2      \sup_{t\leq t_0,\ |x|\geq R}  E^{\lambda^0}\left[
|g(x+X(t))| \right]\\
& =  2\sup_{t\leq t_0,\ |x|\geq R}\Big(  E^{\lambda^0}\left[
|g(x+X(t))| :  |X(t)| < \tfrac{R}{2} \right] +  E^{\lambda^0}\left[%|e^{-S_m^{(j)}(x ,t,;\Phi_m(X))}-e^{-S_0^{(j)}(x,t,;X)} | 
|g(x+X(t))|  :  |X(t)| \geq \tfrac{R}{2} \right]\Big) \\
& \leq 2\left(\sup_{|z|\geq \frac{R}{2}}|g(z)| + \|g\|_\infty \sup_{t \leq t_0}   
\int_{|y|\geq \frac{R}{2}}    k_0^0(y,t)dy   % \lambda^0(|X(t)| \geq \tfrac{R}{2})
 \right).       
\end{align*}
This converges to zero as $R \to \infty$ by  Lemma 5.2, showing claim (i). 

Claim    (ii) can be proved in  the same way as above by using (4.2) and  applying  Lemma 5.1 (ii) and Lemma 5.2.  In fact, we have only to replace $\lambda^0$, $S_{j,A,V}^m(x,t;\Phi_m(X))$, $S_{j,A,V}^0(x,t,X)$, $X(t)$ by  $\mu\times \nu^0$, $S_{3,A,V}^m(x,t; B, T)$, $S_{3,A,V}^0(x,t; B, T)$,         $B(T(t))$,  respectively and note the relation \Equref{probability distribution}.           
$\hfill\square$

\section{Proof of Theorem 1.3}

In this section,  we prove Theorem 1.3.    The proof of $L^2$-convergence for heat semigroups as $m \downarrow 0$ is  not so easy as that of   $C_\infty$-convergence (\cite{I and M 2014}). The reason for this is that,    for example, it is not  trivial that  
\begin{align*}
&\exp\left\{i\int_0^t\int_{0<|z|<1} A\big(x+\Phi_m(X)(s-)+\tfrac{1}{2}\phi_m(z)\big) \cdot \phi_m(z)\ \widetilde{N^0_X}(dsdz)\right\}\\
&  \to  \exp\left\{i\int_0^t\int_{0<|z|<1} A\big(x+X(s-)+\tfrac{1}{2}z\big)\cdot z \ \widetilde{N_X^0}(dsdz)\right\}\quad \text{as }\ m\downarrow 0,
\end{align*} 
since  $A$ may  not be  continuous.   To overcome  this difficulty, we note the following facts: 

(1)  If $0\leq V\in L^1_\text{loc}(\R^d;\R)$, then there exists a sequence $\{V_\ell\}\subset C_0^\infty(\R^d;\R)$ such that 
\begin{align*}
0\leq V_\ell(x)\leq V(x) \text{ a.s.},\ V_\ell \to  V \text{ in } L^1_\text{loc}(\R^d;\R).
\end{align*}

(2)  If $A \in L^{1+\delta}_{\text{loc}}(\R^d;\R^d)$ for some $\delta>0$, then  there exists   a sequence   $\{A_\ell\}\subset C_0^\infty(\R^d;\R^d)$ such that
\begin{align*}
A_{\ell} \to A \ \text {in } L^{1+\delta}_{\text{loc}}(\R^d;\R^d).   
\end{align*}

(3) If $A \in L^{2}_{\text{loc}}(\R^d;\R^d)$, $\nabla \cdot A \in L^{1}_{\text{loc}}(\R^d;\R)$, then  there exists  a sequence   $\{A_\ell\} \subset C_0^\infty(\R^d;\R^d)$ such that
 \begin{align*}
A_{\ell} \to A \ \text {in } L^{2}_{\text{loc}}(\R^d;\R^d),  \   \nabla \cdot A_{\ell} \to \nabla \cdot A \ \text {in } L^{1}_{\text{loc}}(\R^d;\R).
\end{align*}

\begin{lemma} Let  $\{V_\ell\}\subset C_0^\infty(\R^d;\R)$
be an approximate sequence of scalar function $V$ as in (1).  
Then for any  $0<t_0<\infty$, $0<R<\infty$, the following holds: \\ 
 {\rm (i)} Let $\{A_\ell\}\subset C_0^\infty(\R^d;\R^d)$ be an approximate sequence of vector function $A$ as in (2).  Furthermore, let $\{m\}$ be  a decreasing sequence  such that  
 $\displaystyle \sup_{t \leq t_0} |\Phi_m(X)(t)-X(t)| \to 0$ as $m \downarrow 0$ 
$\lambda^0$-a.s.  Then  for $j=1,2$, it holds that as $\ell \to \infty$,  
\begin{align}
  \limsup_{m \downarrow 0}\sup_{t \leq t_0}    \int_{|x|<R}   E^{\lambda^0} \Big[\big|e^{-S^{m}_{j,A,V}(x,t;  \Phi_{m}(X))}- e^{-S^{m}_{j,A_\ell, V_\ell} (x ,t; \Phi_{m}(X))} \big| \Big]dx & \to 0, \\%\quad \text{as } \ \ell \to \infty,     \\    
\sup_{t \leq t_0}    \int_{|x|<R}   E^{\lambda^0} \Big[\big|e^{-S^{0}_{j,A,V}(x,t; X)}- e^{-S^{0}_{j,A_\ell, V_\ell}(x ,t; X)} \big| \Big]dx & \to 0.
\end{align}
 {\rm (ii)} Let $\{A_\ell\}\subset C_0^\infty(\R^d;\R^d)$ be an approximate sequence of vector function $A$ as in (3). 
Then  it holds that as $\ell \to \infty$, 
\begin{align}
  \limsup_{m \downarrow 0}\sup_{t \leq t_0}    \int_{|x|<R}   E^{\mu \times \nu^0} \Big[\big|e^{-S_{3,A,V}(x,t; B,  \Psi_m(T))}- e^{-S_{3,A_\ell, V_\ell}(x ,t; B, \Psi_m(T) )} \big|\Big]dx & \to 0,\\%\quad \text{as } \ \ell \to \infty,     \\    
\sup_{t \leq t_0}    \int_{|x|<R}   E^{\mu \times \nu^0} \Big[\big|e^{-S_{3,A,V}(x,t; B, T)}- e^{-S_{3,A_\ell, V_\ell} (x ,t; B,T)} \big|\Big]dx & \to 0.
\end{align}
\end{lemma}

\noindent
{\itshape Proof.} (i)  
We may  assume without  loss of generality that $0<\delta<1$  because $L_{\text{loc}}^q \subset L_{\text{loc}}^p$ for $1\leq p<q<\infty$.  First, we prove  (6.1)    for $j=1$.  
By \Equref{semimartingale 1},  we have
\begin{align}\label{new estimate}
& \limsup_{m\downarrow 0}\sup_{t \leq t_0} \int_{|x|<R} E^{\lambda^0} \left[|e^{-S_{1,A,V}^m(x,t; \Phi_m(X))}-e^{-S_{1,A_\ell, V_\ell}^{m} (x,t; \Phi_m(X))}|\right] dx \nonumber \\
& \leq  E^{\lambda^0} \left[\limsup_{m\downarrow 0}\sup_{t \leq t_0}\int_{|x|<R}  |e^{-iS_{1,A}^m(x,t; X)}-e^{-iS_{1,A_\ell}^{m} (x,t; X)}| dx\right]  \nonumber \\
& \quad + \limsup_{m\downarrow 0}\sup_{t \leq t_0}  \int_{|x|<R}  E^{\lambda^0}   \left[|e^{-iS_{2,A}^m(x,t;  X)}-e^{-i S_{2,A_\ell}^{m} (x,t;        X)}| \right] dx  \nonumber \\
& \quad +  E^{\lambda^0} \left[\limsup_{m\downarrow 0}\sup_{t \leq t_0}\int_{|x|<R}  |e^{-iS_{3,A}^m(x,t; X)}-e^{-iS_{3,A_\ell }^{m} (x,t; X)}| dx\right]  \nonumber \\
& \quad +  E^{\lambda^0} \left[\limsup_{m\downarrow 0}\sup_{t \leq t_0}\int_{|x|<R}  |e^{-S_{4,V}^m(x,t;  X)}-e^{-S_{4,V_\ell }^{m} (x,t; X)}| dx\right]\nonumber\\
& =: E^{\lambda^0}[I_1^\ell(X)]+\limsup_{m\downarrow 0}\sup_{t\leq t_0}\int_{|x|<R} I_{2}^{m,\ell}(x,t)dx+E^{\lambda^0}[I_3^\ell(X)]+E^{\lambda^0}[I_4^\ell(X)]
\end{align} 
We now show each  term in the last member  of 
\Equref{new estimate}       converges  to zero as $\ell \to \infty$.
To this end, we note that $I_1^\ell(X)$, $I_3^\ell(X)$ and $I_4^\ell(X)$ are less than or equal to $2\mathrm{vol}(R)<\infty$. Here    $\mathrm{vol}(R)$ is the volume of the ball with radius $R$.

{ \it  For the    first   term       of              \Equref{new  estimate}:} 
Since   $|\phi_m(z)|\leq |z|$, we have 
\begin{align*}
I_1^\ell(X)\leq 
&  \sum_{s \leq t_0} \mathbf{1}_{|X(s)-X(s-)| \geq 1}  |X(s)-X(s-)| \int_{|w|<R+ C(X)}|A(w)-A_{\ell}(w)|    dw,
\end{align*}
with a constant $C(X)$ depending on $X$. 
Therefore, here,  since $A_\ell\to A$ in $L^{1+\delta}_{\text{loc}}(\R^d;\R^d)$ and so        
in $L^{1}_{\text{loc}}(\R^d;\R^d)$, it follows that     $E^{\lambda^0}[I_1^\ell(X)]$ converges to zero as $\ell \to \infty$.

{ \it For   the     second   term       of              \Equref{new  estimate}.} 
For convenience of notation, we put 
\begin{align*}
W^{m,\ell}(x,s,y;X):= A(x+\Phi_{m}(X)(s-)+\tfrac{1}{2}y)- A_\ell(x+\Phi_{m}(X)(s-)+\tfrac{1}{2}y). 
\end{align*}
 Let  $G^{ m ,\ell}(x,t;  X)$ be  a subset of   $(0, t] \times \{z; 0< |z| <1\}$     defined by 
\begin{align*}
&G^{ m ,\ell}(x,t;  X):=\Big\{(s,z); \Big| W^{m,\ell}(x,s,\phi_m(z) ;X) \cdot \phi_{m}(z)\Big|>1\Big\}.
\end{align*} 
 Let $\sigma_k(X)$ be the hitting time as    defined in the proof of Lemma 5.1 (i). 
By  the relation $\displaystyle \int_0^t=\int_0^{t\wedge \sigma_k(\Phi_m(X))}+\int_{t \wedge \sigma_k(\Phi_m(X))}^t$,        we have  % $E^{\lambda^0}[\cdots]$ of the second  term on the right-hand side of \Equref{new estimate} is less than or equal to  
\begin{align}\label{estimate second term}
 I_2^{m,\ell}(x,t)& \leq  E^{\lambda^0} \bigg[  \Big|\int_0^{t\wedge \sigma_k(\Phi_m(X))}\int_{0<|z|<1} \mathbf{1}_{G^{ m ,\ell}(x,t;  X)}W^{m,\ell}(x,s,\phi_m(z) ;X) \cdot \phi_m(z) \widetilde{N_X^0} (dsdz) \Big| \bigg]\nonumber\\
& \quad + E^{\lambda^0} \bigg[  \Big|\int_0^{t\wedge \sigma_k(\Phi_m(X))}\int_{0<|z|<1} \mathbf{1}_{G^{ m ,\ell}(x,t;  X)^\complement}W^{m,\ell}(x,s,\phi_m(z) ;X) \cdot \phi_m(z) \widetilde{N_X^0} (dsdz) \Big| \bigg]\nonumber\\
& \quad + 2\lambda^0(\sigma_k(\Phi_m(X))< t) \nonumber\\
&=: J_1^{m,\ell, k }(x,t)+  J_2^{m,\ell, k }(x,t) + 2\lambda^0(\sigma_k(\Phi_m(X))< t).  
\end{align}
For $J_1^{m,\ell, k }(x,t)$, since  $|\widetilde{N_X^0}(dsdz)|\leq N_X(dsdz)+dsn^0(dz)$ and $E^{\lambda^0}\left[ N_X (dsdz)\right]=dsn^0(dz)$, we have  
\begin{align}\label{estimate second term 1} 
  \int_{|x|<R} J_1^{m,\ell, k }(x,t)dx　
& \leq 2 \int_{|x|<R}   E^{\lambda^0} \bigg[ \int_0^{t\wedge \sigma_k(\Phi_m(X))} \int_{0<|z|<1} \mathbf{1}_{G^{ m ,\ell}(x,t;  X)} \nonumber  \\
& \qquad \qquad \qquad \qquad \times \Big|W^{m,\ell}(x,s,\phi_m(z) ;X)  \cdot \phi_{m}(z)  \Big| dsn^{0}(dz)\bigg]dx \nonumber    \\
& \leq 2t_0   \int_{0 < |z| <1} |z|^{1+\delta}  n^{0}(dz)   \int_{|w|< R + k + \frac{1}{2}}|A(w)-A_{\ell}(w)|^{1+ \delta}dw.
\end{align}
For $J_2^{m,\ell, k }(x,t)$, from   the   Schwartz inequality, we have 
\begin{align*}
 J_{2}^{m, \ell, k}(x,t ) ^2        
& \leq  E^{\lambda^0} \bigg[  \int_0^{t\wedge \sigma_k(\Phi_m(X))}\int_{0<|z|<1} \mathbf{1}_{G^{ m ,\ell}(x,t;  X)^\complement}\Big| W^{m,\ell}(x,s,\phi_m(z) ;X) \cdot \phi_m(z) \Big|^{2}      dsn^0(dz)   \bigg]\nonumber\\
& \leq   E^{\lambda^0} \bigg[ \int_{0<|z|<1} |z|^{1+\delta} n^0(dz)  \int_0^{t\wedge \sigma_k(\Phi_m(X))} \Big| W^{m,\ell}(x,s,\phi_m(z) ;X)  \Big|^{1+\delta}      ds  \bigg]. \nonumber
\end{align*}
It follows from  the  Schwartz inequality that 
\begin{align}\label{estimate second term 2}
 \int_{|x|<R} J_{2}^{m,     \ell,k}(x,t) dx
%& \leq \left(|\{|x|<R\}|\int_{|x|<R} J_{2}^{m,     \ell,k}(x,t)^2 dx\right)^{1/2} \nonumber\\
& \leq \bigg(\mathrm{vol}(R)\ t_0 \int_{0 < |z| <1} |z|^{1+\delta}  n^{0}(dz)   \int_{|w|< R + k + \frac{1}{2}}|A(w)-A_{\ell}(w)|^{1+ \delta}dw \bigg)^{1/2}.
\end{align}
By \Equref{estimate second term}, \Equref{estimate second term 1} and \Equref{estimate second term 2}, we have 
\begin{align*}
\limsup_{\ell \to \infty}\limsup_{m\downarrow 0}  \sup_{t \leq t_0}\int_{|x|<R} I_2^{m,\ell}(x,t)  dx & \leq   2
 \lambda^0\left(\limsup_{m\downarrow 0} \Big\{\sigma_k(\Phi_m(X))< t_0\Big\}  \right)    \\
& \leq   2  \lambda^0(\sigma_{k-1} (X)< t_0 ),
\end{align*} 
which converges to zero as $k \to \infty$.   Hence $\displaystyle \limsup_{m\downarrow 0}  \sup_{t \leq t_0}\int_{|x|<R} I_2^{m,\ell}(x,t)  dx$                converges to zero as $\ell \to \infty$.

{ \it For     the     third    term       of              \Equref{new  estimate}:}   Note that  $n^0\phi_m^{-1}(dy)=n^m(dy)=n^m(y)dy$ and  $|\phi_m^{-1}(y)|=\ell_m(|y|)$ (cf. (4.5)). Then
we have  
\begin{align}\label{third term calderon}
\Big|S_{3,A}^m(x,t;  X)-S_{3,A_\ell}^{m} (x,t;        X)\Big|& =  \int_{0}^{t} ds         \text{ p.v.}     \Big|       \int_{0  <|y|<\ell_m^{-1}(1)}W^{m,\ell}(x,s,y ;X) \cdot y \  n^{m}(y) dy       \Big| \nonumber \\
& \leq    \int_{0}^{t} ds          \text{ p.v.}   \int_{0  <|y|<\ell_m^{-1}(1)}\big|  W^{m,\ell}(x,s,y ;X)\big| | y |  \big( n^0(y)-n^m(y)\big)dy     \nonumber \\
& \quad + \int_{0}^{t} ds \text{ p.v.}         \Big|   \int_{0 <|y|<\ell_m^{-1}(1)} W^{m,\ell}(x,s,y ;X) \cdot y \  n^{0}(y) dy                 \Big|     \nonumber \\
&=: K_1^{m, \ell}(x,t;X)+   K_2^{m, \ell} (x,t;X). 
\end{align}
Here we used  the fact that  $n^m(y)< n^0(y)$ (cf. Remark 2.3) in the second inequality. 
For   $K_1^{m,\ell}(x, t;X)$, since        $\int_{|y|>0}(n^{0}(y)-n^{m}(y)) dy=m$ (\cite[Lemma 3.1 (iii)]{I 92}),  we have
\begin{align}\label{third term calderon 1}
&\int_{|x|<R} K_1^{m,\ell}(x,t;X)dx\leq t_0 \ \ell_m^{-1}(1) \   m \int_{|w|<R+C(X)}|A(w)-A_\ell(w)|dw. 
\end{align}
For   $K_2^{m,\ell}(x,      t;X)$,
note that $y=(y_1, \ldots ,y_d) \mapsto y_i n^{0}(y)$  is the Calderon--Zygmund kernel (\cite[p.275]{I 89}) for any $i=1,\ldots ,d$. Then   from  H\"older's  inequality and  the Calderon--Zygmund theorem (\cite[Theorem 2]{S 70}) with a constant $C_{\delta}$ depending only on $\delta$,      we have
\begin{align}\label{third term calderon 2}
&\int_{|x|<R} K_2^{m,\ell}(x,t;X)dx\nonumber \\
& \leq  \mathrm{vol}(R)^{\frac{\delta}{1+\delta}}   \int_0^t ds  \Big(\int_{|x|<R} \text{p.v.} \Big| \int_{0 <|y|<\ell_m^{-1}(1)}   W^{m,\ell}(x,s,y ;X)      \cdot yn^0(y)dy\Big|^{1+\delta}      dx         \Big)^{\frac{1}{1+\delta}}\nonumber \\
&\leq   \mathrm{vol}(R)^{\frac{\delta}{1+\delta}} t_0  C_\delta\Big(\int_{|w|<R+C(X)}|A(w)-A_\ell(w)|^{1+\delta}dw\Big)^{\frac{1}{1+\delta}}. 
\end{align}
By \Equref{third term calderon}, \Equref{third term calderon 1} and \Equref{third term calderon 2},  $I_3^\ell(X)$ converges to zero as  $\ell \to \infty$ and so does $E^{\lambda^0}[I_3^\ell(X)]$.

{\it For   the     fourth    term       of              \Equref{new  estimate}:}  
Since $0\leq V_\ell(x)\leq V(x)$ a.s.,  we have 
 \begin{align*}
&I_4^\ell   (X)  \leq  t_0 \int_{|w| < R+C(X)} (V(w)-V_\ell(w)) dw.
\end{align*} 
It converges to zero as $\ell \to \infty$ since  $V_\ell \to V$ in $L^1_{\text{loc}}(\R^d;\R)$.   Hence $E^{\lambda^0}[I_4^\ell(X)]$      converges to zero as $\ell \to \infty$.

Therefore  we have   (6.1) for $j=1$. Putting $m=0$ in the above proof, we obtain (6.2),  showing claim (i) for $j=1$.

For $j=2$,              (6.1) and  (6.2)  can be proved in the same way as for $j=1$ above. 
This ends  the proof of claim (i).

\vspace{5mm}   
\noindent
(ii) 
By  \Equref{semimartingale 3}, we have 
\begin{align}\label{new estimate 2}
& \limsup_{m\downarrow 0}\sup_{t \leq t_0} \int_{|x|<R} E^{\mu \times \nu^0} \left[|e^{-S_{3,A,V}(x,t; B, \Psi_m(T))}-e^{-S_{3,A_\ell, V_\ell} (x,t;B, \Psi_m(T))}|\right] dx \nonumber \\
&\leq  E^{\nu^0}\Bigg[  \limsup_{m\downarrow 0} \sup_{t \leq t_0} \int_{|x|<R}  E^{\mu}   \left[    |e^{-i S_{1,A}^m(x,t; B, T)}-e^{-iS_{1,A_\ell}^m (x,t;B,  T)}| \right] dx  \Bigg]     \nonumber \\
& \quad +  E^{\mu \times \nu^0} \left[\limsup_{m\downarrow 0}\sup_{t \leq t_0}\int_{|x|<R}  |e^{-iS_{2,A}^m(x,t; B, T)}-e^{-iS_{2,A_\ell }^{m} (x,t;B,  T) }| dx\right]  \nonumber \\
& \quad +  E^{\mu \times \nu^0} \left[\limsup_{m\downarrow 0}\sup_{t \leq t_0}\int_{|x|<R}  |e^{-S_{3,V}^m(x,t; B, T)}-e^{-S_{3,V_\ell }^{m} (x,t;B,  T)}| dx\right]   \nonumber \\
&:= E^{\nu^0}\left[\limsup_{m\downarrow 0}\sup_{t\leq t_0}\int_{|x|<R} I_1^{m, \ell} (x,t; B,T)dx \right]+ E^{\mu\times\nu^0}[I_2^\ell(B,T)]+E^{\mu\times\nu^0}[I_3^\ell(B,T)]. 
\end{align} 
We now   show each term in the last member     of % in the brace $\{ \cdots \}$ of the last member of  
\Equref{new estimate 2}       converges  to zero as $\ell \to \infty$.

{\it For  the    first     term       of              \Equref{new  estimate 2}:} 
By  the  relation   $\displaystyle 
\int_0^t=\int_0^{\Psi_m(T)(t)\wedge \sigma_k(B)}+ \int_{\Psi_m(T)(t)\wedge \sigma_k(B)}^t$ 
  and $(a+b)^2\leq 2a^2+2b^2$ ($a,b \in \R$), we have 
\begin{align*}
 I_1^{m, \ell} (x,t; B,T)^2
 & \leq  2  E^\mu\left[   \int_{0}^{\Psi_m(T) (t) \wedge \sigma_k(B)} | A(x+B(s)) -A_\ell(x+B(s))|^2ds\right] \\
& \quad + 4\mu\big(\sigma_k(B) < \Psi_m(T)(t)\big)^2. 
\end{align*}
It follows  from  the Schwartz inequality  that %,    the integrand of the first term on the right-hand side of \Equref{new estimate 2} is less than or equal to 
\begin{align*}
& \limsup_{m\downarrow 0}\sup_{t\leq t_0}\int_{|x|<R} I_1^{m, \ell} (x,t; B,T)dx\\
& \leq \mathrm{vol}(R)^{1/2} \bigg( 2T(t_0) \int_{|w|<R+k}|A(w)-A_\ell(w)|^2ds+     4\mathrm{vol}(R)   \mu\big(\sigma_k(B)<T(t_0)\big)^2            \bigg)^{1/2}\\
& \to  2 \mathrm{vol}(R) \mu\big(\sigma_k(B)<T(t_0)\big)\quad \text{as }\ \ell \to \infty\\
& \to 0    \quad \text{as }\ k \to \infty.
\end{align*}
 Hence $\displaystyle E^{\nu^0}\left[     \limsup_{m\downarrow 0}\sup_{t\leq t_0}\int_{|x|<R} I_1^{m, \ell} (x,t; B,T)dx\right]$  converges to zero as $\ell \to \infty$.

{\it For   the     second   and third  terms       of              \Equref{new  estimate 2}:} 
Note that $\nabla \cdot A_\ell \to \nabla \cdot A$ and $V_\ell \to V$ in $L^1_\text{loc}(\R^d;\R)$ as $\ell \to \infty$. 
Then we have  
\begin{align*}
& I_2^\ell(B,T) \leq  T(t_0)  \int_{|w| < R+C    (B)} |(\nabla \cdot A)(w)-(\nabla \cdot A_\ell) (w)| dw, \\
&     I_3^\ell(B,T) \leq      t_0 \int_{|w| < R+C(B,T) } \Big(V(w)-V_\ell(w)\Big)  dw,
\end{align*}
 which converge to zero as $\ell \to \infty$. 
 Hence  $E^{\mu\times \nu^0}[I_2^\ell(B,T)]$ and $E^{\mu\times \nu^0}[I_3^\ell(B,T)]$ converge to zero as $\ell \to \infty$.

Therefore we have   (6.3).  Putting $m=0$ in the above proof, we obtain (6.4). 
This ends  the proof of claim  (ii), completing  the proof of Lemma 6.1. 
\hfill$\square$

\vspace{6mm} 

Now we     prove Theorem 1.3.
First, we prove claim (i).  Consider the case $j=1$ or $j=2$. 
Let  $g \in L^2(\R^d)$, $A\in L^{1+\delta}_{\text{loc}}(\R^d;\R^d)$ and    $0\leq V\in L_\text{loc}^1(\R^d;\R)$.  
  Choose a sequence $\{g_n\} \subset C_0^\infty(\R^d)$ such that $g_n \to g$ in $L^2(\R^d)$ as $n \to \infty$.  
Choose sequences  $\{A_\ell \} \subset C_0^\infty(\R^d;\R^d)$ and $\{V_\ell \} \subset C_0^\infty(\R^d; \R)$ as in (2) and (1) at the beginning of this section, respectively.   
 Then we have
\begin{align}\label{last}
& \|e^{-t[H_{j,A}^m-m+V]}g-e^{-t[H_{j,A}^0+V]}g\|_2 \nonumber \\
& \leq \|e^{-t[H_{j,A}^m-m+V]}g-e^{-t[H_{j,A}^m-m+V]}g_n\|_2 
                                                                 + \|e^{-t[H_{j,A}^m-m+V]}g_n-e^{-t[H_{j,A_\ell}^m-m+V_\ell]}g_n\|_2 \nonumber \\
                                                                    & \quad +  \|e^{-t[H_{j,A_\ell}^m-m+V_\ell]}g_n-e^{-t[H_{j,A_\ell}^0+V_\ell]}g_n\|_2 + \|e^{-t[H_{j,A_\ell}^0+V_\ell]}g_n -e^{-t[H_{j,A}^0+V]}g_n\|_2\nonumber \\
                                                                 & \quad + \|e^{-t[H_{j,A}^0+V]}g_n -e^{-t[H_{j,A}^0+V]}g\|_2    \nonumber \\
&=: I_{j}^{m,n}(t) +J_{j}^{m,n, \ell}(t)+K_{j}^{m,n,\ell}(t)+ J_{j}^{0,n,\ell}(t)+I_{j}^{0,n}(t). 
\end{align}
We now estimate each term in the last member  of \Equref{last}.

{\it For  the    first   and  fifth terms       of              \Equref{last}:}   By  the strong continuity  of the semigroup, we have 
\begin{align}\label{i}
  I_{j}^{m,n}(t)+ I_{j}^{0,n}(t)\leq 2\| g_n-g\|_2.%\quad m>0, n\in \mathbf{N}. 
\end{align}
%the first and fifth terms on the right-hand side  of \Equref{last} are less than or equal to    $\|g-g_n\|_2$ by the strong continuity  of the semigroup. 

{\it For  the    third  term      of              \Equref{last}:}  Let $R>0$. From    the  Minkowski     inequality,        we have 
\begin{align}\label{last third term}
K_{j}^{m,n,\ell}(t)
&\leq \|e^{-t[H_{j,A_\ell}^m-m+V_\ell]}g_n-e^{-t[H_{j,A_\ell}^0+V_\ell]}g_n\|_{L^2(|x|<R)} \nonumber \\
& \quad +  \|e^{-t[H_{j,A_\ell}^m-m+V_\ell]}g_n\|_{L^2(|x|\geq R)}+    \|e^{-t[H_{j,A_\ell}^0+V_\ell]}g_n\|_{L^2(|x|\geq R)} \nonumber \\
& \leq \mathrm{vol}(R)^{1/2}  \|e^{-t[H_{j,A_\ell}^m-m+V_\ell]}g_n-e^{-t[H_{j,A_\ell}^0+V_\ell]}g_n\|_{\infty}  \nonumber \\
& \quad + \int_{|x|\geq R} dx\int_{\R^d}     k_0^m(y,t)|g_n(x+y)|^2dy    + \int_{|x|\geq R} dx\int_{\R^d}     k_0^0(y,t)|g_n (x+y)|^2dy. \nonumber 
\end{align}
From Theorem 1.2 (i), the first term in the last member of the above   converges to zero as $m \downarrow 0$ uniformly on $t \leq t_0$.  By the  argument in  \cite[Proof of Theorem 2]{I and M 2014}, the second and third terms in the last member  of the above converges to zero, uniformly on $t \leq t_0$,  $0\leq m \leq 1$. Therefore we have
\begin{align}
\lim_{m\downarrow 0}\sup_{t \leq t_0} K_{j}^{m,n,\ell}(t) =0. 
\end{align}

{\it For  the    second  and fourth  terms      of              \Equref{last}:} Let $R>0$. 
From    the  Minkowski inequality,        we have for $m \geq 0$
\begin{align*}\label{step2-2}
J_{j}^{m,n, \ell}(t) & \leq 
 \|         e^{-t[H_{j,A}^m-m+V]}g_n-e^{-t[H_{j,A_\ell}^m-m+V_\ell]}g_n     \|_{L^2(|x|<R)} \\
& \quad + \|         e^{-t[H_{j,A}^m-m+V]}g_n-e^{-t[H_{j,A_\ell}^m-m+V_\ell]}g_n     \|_{L^2(|x|\geq R)}\nonumber\\
& \leq      \sqrt{2} \|g_n\|_{\infty}  
     \bigg(\int_{|x|<R}   E^{\lambda^0} \Big[\big|e^{-S_{j,A,V}^{m}(x ,t;\Phi_{m}(X))}- e^{-S^{m}_{j,A_\ell, V_\ell}     (x,t; \Phi_{m}(X))} \big|  \Big]dx\bigg)^{1/2}   \nonumber \\
& \quad  + 2 \bigg(\int_{|x| \geq R}dx \int_{\R^d}  k_0^m(y , t) |g_n(x+y)|^2 dy\bigg)^{1/2}. 
\end{align*}
From Lemma 6.1 (i), we have 
\begin{equation}
\begin{aligned}
\limsup_{m\downarrow 0}   \sup_{t\leq t_0} J_{j}^{m,n, \ell}(t) &\to 0 \quad \text{as }\ \ell \to \infty,\\
  \sup_{t\leq t_0} J_{j}^{0,n, \ell}(t) &\to 0 \quad \text{as }\ \ell \to \infty,
\end{aligned}
\end{equation}
where  $\{m\}$ is             a decreasing  sequence such that  $\displaystyle \sup_{t \leq t_0}|\Phi_m(X)(t)-X(t)|\to 0$ as $m \downarrow 0$ $\lambda^0$-a.s.

   Now   let $\{m'\}$ be any subsequence of $\{m\}$ with $m \downarrow 0$. Then,  by Proposition 4.1,  there exists a subsequence $\{m''\}$ of $\{m'\}$ such that $\displaystyle \sup_{t \leq t_0}|\Phi_{m''}(X)(t)-X(t)|\to 0$ as $m \downarrow 0$ $\lambda^0$-a.s. 
By \Equref{last}, (6,14) and  (6.15),  (6.16),  we have 
\begin{align*}
\limsup_{m''\downarrow 0} \sup_{t \leq t_0}\|         e^{-t[H_{j,A}^{m''}-m''+V]}g_n-e^{-t[H_{j,A}^{0}+V]}g_n     \|_{2}
& \leq 2\|g-g_n\|_2, 
\end{align*}
which converges to zero as  $n\to \infty$. This concludes     that        $\displaystyle \sup_{t \leq t_0}\|         e^{-t[H_{j,A}^m-m+V]}g-e^{-t[H_{j,A}^0+V]}g     \|_{2}\to 0$ as $m\downarrow 0$, so   showing claim (i).

Claim (ii) can be proved in the same way as   above by applying Lemma 6.1 (ii), without taking a subsequence $\{m'\}$. 
\hfill$\square$

\vspace{4mm}

\noindent
{\bf Acknowledgment.} 
I would like to thank    Professor  Takashi Ichinose and Professor  Hidekazu   Ito   for many valuable comments,  helpful discussions  and  warm encouragements    during this work.

\vspace{5mm} 

\noindent 
Taro Murayama\\
%Kanazawa University\\
Institute of Science and Engineering,\\
 Kanazawa University,\\
 Kakumamachi, Kanazawa, Ishikawa, 920--1192, Japan\\

\noindent
Current  Address:
Department of General Education, \\
National Institute of Technology, Ishikawa College, \\
Kitachujo, Tsubata, Ishikawa, 929--0392, Japan

\end{document}